\documentclass[12pt,twoside]{article}
\date{}
\usepackage{amsthm, amsmath, amssymb, url}

\usepackage{amssymb,amsfonts}

\def\RR{\mathbb{R}}    
\def\NN{\mathbb{N}}    
\def\ZZ{\mathbb{Z}}    
\def\PP{\mathbb{P}}    
\def\EE{\mathbb{E}}    

\def\proof{\noindent{\em Proof.}~}
\def\endproof{\hfill\mbox{\fbox{\rule{0mm}{0.0mm}}}\vspace{2ex}}
\newtheorem{theorem}{Theorem}
\newtheorem{lemma}[theorem]{Lemma}
\newtheorem{proposition}[theorem]{Proposition}
\newtheorem{corollary}[theorem]{Corollary}
\title{On the rate of convergence of Krasnosel'ski\v{\i}-Mann iterations and their connection with sums of Bernoullis}
\author{R. Cominetti\thanks{Departamento Ingenier\'ia
    Industrial, Universidad de Chile. e-mail: rccc@dii.uchile.cl.
    Supported by Fondecyt 1100046 and Nucleo Milenio Informaci\'on y
    Coordinaci\'on en Redes ICM/FIC P10-024F.} \and 
J.A. Soto\thanks{Departamento Ingenier\'ia Matem\'atica and Centro de Modelamiento Matem\'atico (UMI 2807 CNRS), Universidad de Chile. e-mail: jsoto@dim.uchile.cl. Supported by Basal-Conicyt project and Nucleo Milenio Informaci\'on y Coordinaci\'on en Redes ICM/FIC P10-024F.}
\and J. Vaisman\thanks{Departamento de Ingenier\'\i a Matem\'atica, Universidad de Chile.
e-mail: hellovaisman@gmail.com.}
}
\begin{document}
\maketitle

\begin{abstract}
In this paper we establish an estimate for the rate of convergence of the
Krasnosel'ski\v{\i}-Mann iteration for computing fixed points of 
non-expansive maps. Our main result settles the Baillon-Bruck conjecture \cite{bb}
on the asymptotic regularity of this iteration. The proof proceeds by establishing a 
connection between these iterates and a stochastic process involving  sums of
non-homogeneous Bernoulli trials. We also exploit a new Hoeffding-type
inequality to majorize the expected value of a convex function of these sums 
using Poisson distributions.
\end{abstract}

\vspace{2ex}
\begin{small}
  \\[2ex]
  \noindent{\bf Math Subject Class MSC2010:} 47H09, 47H10, 65K15, 33CXX, 60E15, 60G50
\end{small}

\vfill

\pagebreak
\section{Introduction}

Let $T:C\to C$ be a non-expansive map defined on a  convex subset $C\subseteq X$ 
of a normed space $(X,\|\cdot\|)$. The Krasnosel'ski\v{\i}-Mann iteration for 
computing a fixed-point of $T$ is defined by ({\em cf.} \cite{kra, man}) 
\begin{equation}\label{mann}
x_k=(1-\alpha_k)x_{k-1}+\alpha_k Tx_{k-1}
\end{equation}
with $x_0\in C$ given and $\alpha_k\in [0,1]$. 

Strong convergence of $x_k$ to a fixed point was proved in \cite[Krasnosel'ski\v{\i}]{kra} for 
$\alpha_k\equiv\frac{1}{2}$, when $X$ is a uniformly convex Banach space and $T(C)$ is contained in a compact subset of $C$. 
This result was extended to $\alpha_k\equiv\alpha$ \cite[Schaefer]{sch}  and $X$ strictly convex
\cite[Edelstein]{ede}, while \cite[Ishikawa]{ish} proved it for
general Banach spaces with $\alpha_k$  bounded away from 1
and $\sum\alpha_k\!=\!\infty$. The Banach case with $\alpha_k\equiv\alpha$ was also considered in \cite[Edelstein and O'Brien]{eob}.
Without the compactness assumption, weak convergence
was established in \cite[Reich]{rei} assuming $\sum\alpha_k(1\!-\!\alpha_k)=\infty$ and
$\mbox{\rm Fix}(T)\neq\phi$, for $X$ uniformly convex with a Fr\'echet
differentiable norm. Although strong convergence does not hold in general
(see \cite[Genel and Lindenstrauss]{gen} and \cite[Bauschke {\em et al.}]{bay}), it does occur for most 
operators in the sense of Baire's categories (see \cite[Reich and Zaslavski]{res}).

The crucial step in proving the convergence of the iterates in all these results is to show that $\|x_n-Tx_n\|$ 
tends to 0, a property which is now known as {\em asymptotic regularity} \cite{bbr,bor,bro,rsh}. 
Under various assumptions, asymptotic regularity was also proved in 
\cite[Groetsch]{gro} and \cite[Goebel and Kirk]{goe}. The latter noted a certain uniformity in the convergence, 
namely, for each $\epsilon>0$ we have $\|x_n-Tx_n\|\leq \epsilon$ for all $n\geq n_0$, with $n_0$ 
depending on $\epsilon$ and $C$ but independent of the initial point $x_0$ and the map $T$.
More recently, using proof mining techniques, Kohlenbach \cite{koh,ko1} showed that $n_0$ could be 
chosen to depend on $C$ only through its diameter.
An explicit metric estimate which readily implies all
these results was stated in \cite[Baillon and Bruck]{bb},
namely, they conjectured the existence of a universal constant $\kappa$ such that 
\begin{equation}\label{bnd}
\|x_n-Tx_n\|\leq \kappa\frac{\mbox{diam}(C)}{\sqrt{\sum_{i=1}^n\alpha_i(1\!-\!\alpha_i)}}
\end{equation}
and proved it for the case $\alpha_i\equiv\alpha$  with $\kappa=\frac{1}{\sqrt{\pi}}$.

In this paper we settle this conjecture by proving that the bound holds in general 
with $\kappa=\frac{1}{\sqrt{\pi}}$ for any sequence $\alpha_k$ and each  non-expansive $T:C\to C$. 
Although we do not know whether this is the smallest possible  $\kappa$, we provide an example which shows that it
cannot be improved 
by more than 17\%. We also discuss how the result
can be used to analyze the convergence of (\ref{mann}), and 
how it applies when $C$ is unbounded but $\mbox{Fix}(T)\neq\phi$.

Our proof  is based on a recursive bound 
for the distances between the iterates $\|x_m-x_n\|\leq c_{mn}$, 
where $c_{mn}$ admits a nice probabilistic interpretation
in terms of a random walk on $\ZZ$. In proving the theorem we 
exploit some properties of 
the hypergeometric and modified Bessel functions, as well as 
a known identity for Catalan numbers. We also use the following Hoeffding-type 
inequality  which might be of interest on its own: 
{\em if $S=X_1+\cdots+X_m$ is a sum of independent Bernoullis and $Z$ is a 
Poisson with the same mean $\EE(Z)=\EE(S)$, then $\EE[g(S)]\leq\EE[g(Z)]$
for every convex function $g:\NN\to\RR$.}

\section{Main result}\label{smain}

\begin{theorem}\label{main}
The Krasnosel'ski\v{\i}-Mann iterates generated by {\em (\ref{mann})} satisfy
\begin{equation}\label{bnd0}
\|x_n-Tx_n\|\leq \frac{\mbox{\rm diam}(C)}{\sqrt{\pi\sum_{i=1}^n\alpha_i(1\!-\!\alpha_i)}}.
\end{equation}
\end{theorem}

The proof is split into several intermediate steps. 
Note that by rescaling the norm, we may assume $\mbox{diam}(C)\!=\!1$.
\vspace{-1ex}

\subsection{A recursive bound}
\vspace{-0.5ex}
Let $\rho_k=\Pi_{j=1}^k(1\!-\!\alpha_j)$  and
$\pi_k^n=\rho_n\frac{\alpha_k}{\rho_k}=\alpha_k\Pi_{j=k+1}^n(1\!-\!\alpha_j)$. 
By convention we also set $\rho_0\!=\!\alpha_0\!=\!1$, while the term $Tx_{-1}$ is interpreted as $x_0$.
\begin{proposition}
For $n\geq 0$ we have $x_n=\mbox{$\sum_{k=0}^n\pi_k^nTx_{k-1}$}$
	and \vspace{-0.5ex}
\begin{equation}\label{ddd}
x_m-x_n=
\sum_{j=0}^m\sum_{k=m+1}^n\!\!\!\mbox{$\pi_j^m\pi_k^n$}\,[Tx_{j-1}-Tx_{k-1}]\quad  \mbox{for $0\leq m\leq n$.}
\end{equation}
\end{proposition}
\proof
Dividing (\ref{mann})  by $\rho_k$ we have
$\mbox{$\frac{x_k}{\rho_k}=\frac{x_{k-1}}{\rho_{k-1}}+\frac{\alpha_k}{\rho_k}Tx_{k-1}$}$
which, when iterated, yields
$\mbox{$\frac{x_n}{\rho_n}$} =x_0+\mbox{$\sum_{k=1}^n\frac{\alpha_k}{\rho_k}Tx_{k-1}$}.$
Using the conventions $\rho_0=\alpha_0=1$ and $x_0=Tx_{-1}$ we get precisely
$x_n=\sum_{k=0}^n\pi^n_kTx_{k-1}$. 
This equality, combined with the identities $\sum_{j=0}^m\pi_j^m=1$ and $\pi_k^m-\pi_k^n=\sum_{j=m+1}^n\pi_j^n\pi_k^m$,  
yields\vspace{-0.5ex}
\begin{eqnarray*}
x_m-x_n&=&\sum_{k=0}^m(\pi_k^m\!-\!\pi_k^n)Tx_{k-1}-\sum_{k=m+1}^n\pi_k^nTx_{k-1}\\
&=&\sum_{k=0}^m\sum_{j=m+1}^n\pi_j^n\pi_k^mTx_{k-1}
-\sum_{j=0}^m\sum_{k=m+1}^n\pi_j^m\pi_k^nTx_{k-1}
\end{eqnarray*}
so that exchanging $j$ and $k$ in the first double sum 
we obtain (\ref{ddd}).
\endproof

\begin{corollary} \label{coro}
Define  $c_{mn}$ recursively by setting $c_{-1,n}=1$  for all $n\geq 0$ and
$$c_{mn}\hspace{1ex}=\mbox{$\sum_{j=0}^m\sum_{k=m+1}^n$}\pi_j^m\, \pi_k^n\, c_{j-1,k-1}\quad\mbox{for $0\leq m\leq n$}.\leqno(R)
$$
Then $\|x_m-x_n\|\leq c_{mn}$ for all $0\leq m\leq n$.
\end{corollary}
\proof The proof is by induction on $n$. Suppose that $\|x_j-x_k\|\leq c_{jk}$ holds for all
$0\leq j\leq k\leq n-1$. Using the triangle inequality in  (\ref{ddd}) we get
\begin{equation}\label{cota}
\|x_m-x_n\|\leq\mbox{$\sum_{j=0}^m\sum_{k=m+1}^n$}\pi_j^m \pi_k^n \|Tx_{j-1}-Tx_{k-1}\|.
\end{equation}
The induction hypothesis gives $\|Tx_{j-1}-Tx_{k-1}\|\leq  \|x_{j-1}-x_{k-1}\|\leq c_{j-1,k-1}$ for $1\leq j < k$,
while for $j=0$ we have
$\|Tx_{-1}-Tx_{k-1}\|=\|x_0-Tx_{k-1}\|\leq\mbox{diam}(C)=1=c_{-1,k-1}.$
Plugging these bounds into (\ref{cota}) and using $(R)$ we
deduce $\|x_m-x_n\|\leq c_{mn}$ completing the induction step.
\endproof

\vspace{2ex}
Note that for $m=n$ we have $c_{nn}=0$ and the inequality $\|x_n-x_n\|\leq c_{nn}$ holds trivially.  
More interestingly, since $\|x_n-x_{n+1}\|=\alpha_{n+1}\|x_n-Tx_n\|$ we have
$\|x_n-Tx_n\|\leq \mbox{$\frac{c_{n,n+1}}{\alpha_{n+1}}$} \triangleq P^n$
so that
Theorem \ref{main} will follow by showing
\begin{equation}\label{bnd00}
\mbox{$\sqrt{\sum_{i=1}^n\alpha_i(1\!-\!\alpha_i)}\;P^n\leq\frac{1}{\sqrt\pi}$}.
\end{equation}
Our analysis proves 
that this bound is sharp, so that $\frac{1}{\sqrt{\pi}}$ is the
best constant one can get from Corollary \ref{coro}.
This does not exclude the possibility that other techniques might lead to sharper bounds
in Theorem \ref{main} ({\em cf.} \cite[Baillon and Bruck]{bb2}).

\subsection{Fox-and-Hare race and a random walk}

The recurrence $(R)$ has a
probabilistic interpretation. Consider a fox at
position $n$ trying to catch a hare located at $m<n$.
At each integer $i\in \NN$ the fox must jump over a hurdle to reach 
$i\!-\!1$. The jump succeeds with probability $(1\!-\!\alpha_i)$ in which 
case the process repeats, otherwise the fox falls at $i\!-\!1$ where it
rests to recover from injuries. Thus, starting from $n$ the probability of landing at $k\!-\!1$ is 
precisely $\pi_k^n$. The fox catches the hare if it jumps successfully down to $m$ or below. 
Otherwise,  the hare runs toward the burrow located at $-1$ by
following the same rules. The process alternates  until either the fox catches the hare,
or the hare reaches the burrow. 

The recurrence $(R)$ satisfied by $c_{mn}$ characterizes precisely 
the probability for the hare to reach the burrow safely when the process starts at $(m,n)$. 
This is also consistent with the boundary cases
$c_{-1,n}=1$ and $c_{nn}=0$. Note that $\alpha_0=1$ so at $i=0$ 
both the fox and hare fall with certainty, landing at $-1$.
From this interpretation we get the following 
expression for $c_{mn}$.
\begin{proposition} \label{bbb}
Let $(F_i)_{i\in\NN}$ and $(H_i)_{i\in\NN}$ denote independent Bernoulli trials representing respectively the events
that the fox and hare fail at the $i$-th hurdle, so that 
$\PP(F_i\!=\!1)\!=\!\PP(H_i\!=\!1)\!=\!\alpha_i$. Then
\begin{equation}\label{jumps}
c_{mn}=\PP(\mbox{$\sum_{i=k}^nF_i>\sum_{i=k}^mH_i$ for all $k=m+1,\ldots,1$}).
\end{equation}
In particular, denoting $Z_i=F_i\!-\!H_i$ we have
\begin{equation}\label{jumps2}
P^n=\mbox{$\frac{c_{n,n+1}}{\alpha_{n+1}}$}=\PP(\mbox{$\sum_{i=k}^nZ_i\geq 0$ for $k=n,\ldots,1$}).
\end{equation}
\end{proposition}
\proof 
Formula (\ref{jumps}) is just a restatement of the fact that the hare wins iff
the number of times the fox falls in any interval $\{k,\ldots,n\}$
is strictly larger than the number of falls of the hare in $\{k,\ldots,m\}$.
The expression for $P^n$ follows by noting that
the event corresponding  to $c_{n,n+1}$ in (\ref{jumps}) requires $F_{n+1}=1$ (take $k=n+1$).
\endproof

Formula (\ref{jumps2}) has an alternative interpretation. Let  $p_i\!=\!2\alpha_i(1\!-\!\alpha_i)$
so that $Z_i$ takes values in $\{-1,0,1\}$ with probabilities $\frac{p_i}{2},1\!-\!p_i,\frac{p_i}{2}$.
The sums $\sum_{i=k}^nZ_i$ taken in reverse order $k=n,\ldots,1$
define a random walk on $\ZZ$ where at each stage the process 
stays at the current  position with some probability, and otherwise 
moves left or right with equal probability  
as in a standard random walk. Hence, $P^n$ is the probability that
the walk remains non-negative over $n$ stages. Conditioning on the total number of stages at which the process effectively moves,
this is also the probability that a standard random walk stays non-negative over a 
random number of stages. Using this interpretation we get the following  more explicit  formula.

\begin{proposition} \label{sb}
Let $M=M_1+\ldots+M_n$ be a sum of independent Bernoullis
with success probabilities $\PP(M_i\!=\!1)\!=\!p_i\!=\!2\alpha_i(1\!-\!\alpha_i)$ and consider the integer 
function 
$F(m)= {m \choose \lfloor m/2\rfloor}2^{-m}$.
Then $P^n=\EE[F(M)]$.
\end{proposition}
\proof 
The variable $M_i$ can be interpreted as {\em move/stay} and $Z_i$ can be expressed as $Z_i=M_iD_i$ with 
$D_i$ independent variables representing the direction of the movement: $\PP(D_i\!=\!-1)=\PP(D_i\!=\!1)=\frac{1}{2}$.
Conditioning on the sum $M$ and using the exchangeability of the variables $D_i$ 
we obtain
\begin{eqnarray*}
P^n&=&\mbox{$\sum_{m=0}^n$}\PP(\mbox{$\sum_{i=k}^nM_iD_i\geq 0$ for $k=n,\ldots,1$}|M=m)\PP(M=m)\\
&=&\mbox{$\sum_{m=0}^n$}\PP(\mbox{$\sum_{j=1}^\ell D_j\geq 0$ for $\ell=1,\ldots,m$})\PP(M=m).
\end{eqnarray*}
The expression $\PP(\mbox{$\sum_{j=1}^\ell D_j\geq 0$ for $\ell=1,\ldots,m$})$ is the probability that 
a standard random walk started from 0 remains non-negative over $m$ stages. Its value is
precisely $F(m)$ \cite[Ch. III.3]{fel} so the conclusion follows.
\endproof

The next result establishes an alternative recursion satisfied by $c_{mn}$.
This is not used in our proof, but we state
in case someone could use it to find a simpler proof of Theorem \ref{main}.
\begin{proposition} Denoting $\bar\alpha_k=1\!-\!\alpha_k$, we have the recurrence
\begin{equation}\label{recurrence}
c_{mn}=\bar\alpha_mc_{m-1,n}+\bar\alpha_n c_{m,n-1}+(\alpha_n\alpha_m-\bar\alpha_n\bar\alpha_m)c_{m-1,n-1}.
\end{equation}
\end{proposition}
\proof Denote $w_{jk}=\pi_j^m\pi_k^nc_{j-1,k-1}$ and let $S=A+B-C-D$ with
$$\begin{array}{cclcl}
A&=&c_{mn}&=&\mbox{$\sum_{j=0}^{m}\sum_{k=m+1}^nw_{jk}$}\\
B&=&\bar\alpha_m\bar\alpha_nc_{m-1,n-1}&=&\mbox{$\sum_{j=0}^{m-1}\sum_{k=m}^{n-1}w_{jk}$}\\
C&=&\bar\alpha_mc_{m-1,n}&=&\mbox{$\sum_{j=0}^{m-1}\sum_{k=m}^nw_{jk}$}\\
D&=&\bar\alpha_nc_{m,n-1}&=&\mbox{$\sum_{j=0}^{m}\sum_{k=m+1}^{n-1}w_{jk}$}.
\end{array}
$$
Canceling out the common terms we get
$S=w_{mn}=\alpha_m\alpha_nc_{m-1,n-1}$
which is exactly (\ref{recurrence}).
\endproof

\subsection{A sharp upper bound}\label{dostres}
From Proposition \ref{sb}, the bound (\ref{bnd00}) is equivalent to showing that
$$R^n(p)\triangleq\sqrt{p_1+\ldots+p_n}\; \EE[F(M_1+\ldots+M_n)]\leq \mbox{\scriptsize $\sqrt{\frac{2}{\pi}}$}$$
for all $n$ and $0\leq p_i\leq \frac{1}{2}$.
The function $R^n(p)$ is strictly concave in each 
variable $p_i$ separately, so the maximum is attained at the 
extreme values $0,\frac{1}{2}$ or at a unique point in $(0,\frac{1}{2})$. Interestingly, 
all non-extreme coordinates may be taken equal.
\begin{lemma}\label{lemita}
$R^n(p)$ is maximal when $p_i\in\{0,u,\frac{1}{2}\}$ for some $0<u<\frac{1}{2}$.
\end{lemma}
\proof
Let $p$ maximize $R^n(p)$ and suppose $p_j=x$  and $p_k=y$ with $x,y\in (0,\frac{1}{2})$
and $x\neq y$.
Let $h(k)=\EE[F(k+S)]$ where $S=\sum_{i\neq j,k}M_i$ so that
\begin{eqnarray*}
P^n&=&(1\!-\!x)(1\!-\!y)h(0)+[x(1\!-\!y)+y(1\!-\!x)]h(1)+xyh(2)\\
&=&a+b(x\!+\!y)+cxy
\end{eqnarray*}
with $a\!=\!h(0)$, $b\!=\!h(1)\!-\!h(0)$ and $c\!=\!h(0)\!+\!h(2)\!-\!2h(1)$.
Setting $m\!=\!\!\sum_{i\neq j,k}p_i$ it follows that $x,y\in (0,\frac{1}{2})$ 
maximize  the expression 
$$\sqrt{m+x+y}\;[a+b(x+y)+cxy].$$
Setting the partial derivatives to 0 we get $cx=cy$
and since $x\neq y$ it follows that $c=0$. But then, the function 
depends only on the sum $x+y$ and we may change these 
coordinates to $x+\epsilon, y-\epsilon$ keeping the same value, 
until one of them hits an extreme value: either $x+\epsilon=\frac{1}{2}$
or $y-\epsilon=0$. This yields a new optimal $p$ with
one coordinate less in $(0,\frac{1}{2})$.  Repeating this process 
we get an optimal $p$ whose coordinates take at most one 
value in $(0,\frac{1}{2})$.
\endproof


According to this Lemma, in order to bound $R^n(p)$ it suffices to consider
the case $p_i\in \{0,u,\frac{1}{2}\}$ with $0<u<\frac{1}{2}$.
Moreover, by changing $n$ we may ignore the deterministic variables with $p_i=0$. 
We distinguish two cases.

\subsubsection{All coordinates $p_i=u$}
In this case  $R^n(p)=\sqrt{nu}\;\EE[F(S)]$ with $S\sim B(n,u)$ Binomial. 
This case follows from the results in \cite[Baillon and Bruck]{bb}
which were obtained using a computer generated proof.
Here we provide a direct proof based on a known identity for Catalan numbers.
\begin{proposition} \label{equal} Let $S\sim B(n,u)$ with $0<u<\frac{1}{2}$. Then
\begin{equation}\label{beq}
\EE[F(S)]=\mbox{$\sum_{k=0}^{n}\frac{(-1)^k}{k+1}{2k\choose k}{n \choose k}(\frac{u}{2})^k$}
\end{equation}
and $R^n(p)=\sqrt{nu}\;\EE[F(S)]$ increases with $n$ towards \mbox{\scriptsize $\sqrt{\frac{2}{\pi}}$}.
\end{proposition}
\proof
Using the Binomial theorem, a straightforward computation gives 
\begin{eqnarray}\nonumber
\EE[F(S)]&=&\mbox{$\sum_{j=0}^nF(j){n \choose j}u^j(1\!-\!u)^{n-j}$}\\ \nonumber
&=&\mbox{$\sum_{j=0}^nF(j){n \choose j}u^j\sum_{i=0}^{n-j}{n-j\choose i}(-u)^i$}\\ \nonumber
&=&\mbox{$\sum_{j=0}^n\sum_{k=j}^n(-1)^jF(j){n \choose j}{n-j\choose k-j}(-u)^k$}\\ \label{54}
&=&\mbox{$\sum_{k=0}^{n}{n \choose k}(-u)^k\sum_{j=0}^k(-1)^j{k\choose j}F(j)$}
\end{eqnarray}
where the last equality follows from the 
identity ${n\choose j}{n-j\choose k-j}={n\choose k}{k\choose j}$ and exchanging the order of the sums.
The last inner sum  may be computed from a known identity for Catalan numbers 
$C_k=\frac{1}{k+1}{2k\choose k}$, namely\footnote{See
{\scriptsize \tt http://mathworld.wolfram.com/CatalanNumber.html}. A proof is also given in \S\ref{acatalan}.}
$$C_k=\mbox{$\sum_{j=0}^k(-1)^j2^{k-j}{k\choose j}{j\choose \lfloor j/2\rfloor}$}
=2^k\mbox{$\sum_{j=0}^k(-1)^j{k\choose j}F(j)$}$$
which when substituted into (\ref{54}) yields (\ref{beq}).

By direct verification, the
 expression on the right of (\ref{beq})  is the hypergeometric function $_2F_1(-n,\frac{1}{2};2;2u)$, whose
Euler integral representation gives
$$\EE[F(S)]=\mbox{$\frac{2}{\pi}\int_0^1 t^{-1/2}(1\!-\!t)^{1/2}(1\!-\!2ut)^ndt$}.$$
Multiplying by $\sqrt{nu}$ and using the change of variables $s=2nut$ we get
$$R^n(p)=\sqrt{nu}\;\EE[F(S)]=\mbox{$\frac{\sqrt{2}}{\pi}$}\int_0^{2nu}\!\!\!
\mbox{$\sqrt{\frac{1}{s}-\frac{1}{2nu}}$}\;\mbox{$(1\!-\!\frac{s}{n})^nds$}$$
which increases with $n$ towards the limit 
$\frac{\sqrt{2}}{\pi}\int_0^\infty\!\!\frac{1}{\sqrt{s}}e^{-s}ds=\frac{\sqrt{2}}{\pi}\,\Gamma(\frac{1}{2})=\mbox{\scriptsize $\sqrt{\frac{2}{\pi}}$}$.
\endproof

\subsubsection{At least one coordinate $p_i=\frac{1}{2}$}
With no loss of generality assume $p_1=\frac{1}{2}$ and denote $S=M_2+\ldots+M_{n}$.
Conditioning on $M_1$ and setting $g(k)\triangleq \mbox{$\frac{1}{2}$}[F(k)+F(k\!+\!1)]$ we get
$$\EE[F(M_1+\ldots+M_n)]=\EE[g(S)].$$
A direct calculation shows that $g:\NN\to\RR$ is convex, namely
$$g(k)\leq\mbox{$\frac{1}{2}$}[g(k\!-\!1)+g(k\!+\!1)]\quad\mbox{for all $k\geq 1$},$$
so we may use the Hoeffding-type inequality in Proposition \ref{poisson}
to obtain $\EE[g(S)]\leq\EE[g(Z)]$ with $Z\sim P(z)$ a Poisson variable  with $z=p_2+\cdots+p_n$.
From this it follows that
 \begin{eqnarray}\nonumber
R^n(p)&\leq&\mbox{$\sqrt{z\!+\!\frac{1}{2}}$}\;\EE[g(Z)]\\
&=&\mbox{$\frac{1}{2}\sqrt{z\!+\!\frac{1}{2}}$}\mbox{$\;\sum_{k=0}^\infty [F(k)+F(1+k)]\exp(-z)\frac{z^k}{k!}.$} \nonumber\\
&=&\mbox{$\sqrt{z\!+\!\frac{1}{2}}$}\;\exp(-z)\mbox{$[I_0(z)+(1-\frac{1}{2z})I_1(z)]$}\label{cotafinal}
\end{eqnarray}
where
$I_0(z)=\mbox{$\sum_{k=0}^\infty\frac{1}{(k!)^2}(\frac{z}{2})^{2k}$}$
and $I_1(z)=\mbox{$\sum_{k=0}^\infty\frac{1}{k!(k+1)!}(\frac{z}{2})^{2k+1}$}$
are modified Bessel functions.

\begin{proposition} Let $h(z)$ denote the expression in {\em (\ref{cotafinal})}.
Then $h(z)$ is increasing with $h(z)\leq \lim_{z\to\infty}h(z)=\mbox{\scriptsize $\sqrt{\frac{2}{\pi}}$}$. 
\end{proposition}
\proof The identities $I_0'(z)=I_1(z)$ and $I_1'(z)=I_0(z)-\frac{1}{z}I_1(z)$ imply
$$h'(z)=\mbox{$\frac{\exp(-z)}{4z^2\sqrt{z+\frac{1}{2}}}[2(1+z)I_1(z)-zI_0(z)]$}$$
so that proving that $h$ is increasing reduces to $zI_0(z)\leq 2(1+z)I_1(z)$.
Letting $x=z/2$ and rearranging terms, this is equivalent to
$$\mbox{$\sum_{k=1}^\infty\frac{x^{2k+1}}{(k-1)!(k+1)!}$}
\leq \mbox{$2\sum_{k=0}^\infty\frac{x^{2k+2}}{k!(k+1)!}$}.$$
This latter inequality follows easily by noting that each term on the left can be
bounded from above by two consecutive terms on the right, namely
$$\mbox{$\frac{x^{2k+1}}{(k-1)!(k+1)!}\leq \frac{x^{2k}}{(k-1)!k!}+\frac{x^{2k+2}}{k!(k+1)!}$}$$
which results from the trivial inequality $kx\leq k(k+1)+x^2$.

Thus $h(z)$ is increasing and therefore it is bounded from
above by its limit $\ell=\lim_{z\to\infty}h(z)$. To prove that $\ell=\mbox{\scriptsize $\sqrt{\frac{2}{\pi}}$}$
one may use the known asymptotics
$\exp(-z)\sqrt{z}\;I_\alpha(z)\to \mbox{\scriptsize $\frac{1}{\sqrt{2\pi}}$}$  (see \cite[Chapter 9]{abr}).
Alternatively, one may use the integral representation 
$I_n(z)=\frac{1}{\pi}\int_0^\pi \cos(n\theta)e^{z\cos\theta} d\theta$
to write
$$\ell=\lim_{z\to\infty}\mbox{$\frac{1}{\pi}$}
\mbox{$\sqrt{z\!+\!\frac{1}{2}}$}\int_0^\pi\!\! \mbox{$[1+(1\!-\!\frac{1}{2z})\cos\theta]e^{-z(1-\cos\theta)}d\theta$}.$$
Since $\frac{1}{2z}\mbox{\scriptsize $\sqrt{z\!+\!\frac{1}{2}}$}\to 0$ the relevant term for the limit  
is $\int_0^\pi [1+\cos\theta]e^{-z(1-\cos\theta)}d\theta$, which is transformed 
by the change of variables $z(1\!-\!\cos\theta)=x^2/2$ into
$$\ell=\lim_{z\to\infty}\mbox{$\frac{2}{\pi}\mbox{$\sqrt{1\!+\!\frac{1}{2z}}$}$}\int_0^{\sqrt{4z}}\!\!\!\!\!\!
\mbox{$(1\!-\!\frac{x^2}{4z})^{1/2}$}\,e^{-x^2/2}\;dx
=\mbox{$\frac{2}{\pi}$}\!\int_0^{\infty} \!\!\!\!e^{-x^2/2} dx = \mbox{\scriptsize $\sqrt{\frac{2}{\pi}}$}.\vspace{-4ex}
$$
\endproof

\vspace{2ex}
\noindent{\sc Remark.} An alternative proof of the monotonicity of $h(z)$
is obtained by substituting the well-known recurrence $I_{n+1}\!=\!I_{n-1}\!-\!\frac{2n}{z}I_n$
into the Turan-type inequality $I_{n-1}I_{n+1}\leq I_n^2$  (see \cite{thi}) which gives
$I_{n-1}^2 - \mbox{$\frac{2n}{z}$} I_{n-1} I_n \leq I_n^2$. Denoting $x=I_{n-1}/I_n$ we have
$x^2-\frac{2n}{z}x\leq 1$, and solving the quadratic we get
 $\mbox{$ x \leq \frac{n}{z} + \sqrt{ 1 + (\frac{n}{z})^2}.$}$
For $n=1$ this last expression is smaller than $2(z+1)/z$
which gives  $z I_0(z) \leq  2 (z+1) I_1(z)$ so that
$h'(z)\geq 0$.

\subsection{Conclusion}
The bounds in \S\ref{dostres} establish (\ref{bnd00}) and prove Theorem \ref{main}.
Moreover,  the bound (\ref{bnd00}) is sharp and cannot be improved.
Indeed, for $\alpha_i\equiv\alpha$ constant, 
setting $u=2\alpha(1\!-\!\alpha)$ and $S\sim B(n,u)$ we have
$$ \mbox{$\sqrt{\sum_{i=1}^n\alpha_i(1\!-\!\alpha_i)}\;P^n=
\mbox{$\sqrt{\frac{nu}{2}}$}\;\EE[F(S)]$}$$
and by Proposition \ref{equal}  this quantity converges
to \mbox{\scriptsize $\frac{1}{\sqrt{\pi}}$} as $n\to\infty$.
This does not mean that (\ref{bnd0}) is itself sharp
since we only have $\|x_n-Tx_n\|\leq P^n$. Thus, a natural question 
is to find the smallest constant $\kappa$ for which (\ref{bnd}) holds.
Although we do not know whether (\ref{bnd0}) is 
sharp or not, the following example shows that this bound cannot be
improved by more than 17\%.

\vspace{2ex}
\noindent{\sc Example.} Take $X=\ell^1(\NN)$ and let $C$ be the set of all 
sequences $x=(x^i)_{i\in\NN}$ with $x^i\geq 0$ and $\sum_{i=0}^\infty x^i\leq 1$,
 so that diam$(C)=2$.
Let $T:C\to C$ be the right-shift isometry 
$T(x^0,x^1,x^2,\ldots)=(0,x^0,x^1,x^2,\ldots)$. 
Then, the iteration $(KM)$ started from $x_0=(1,0,0,\ldots)$ generates a sequence of the form $x_n=(p^0_n,p^1_n,\ldots,p_n^n,0,0,\ldots)$ with
$$p^i_n=\PP(X_1+\ldots+X_n=i)$$
where $X_i$ are independent Bernoullis with $\PP(X_i=1)=\alpha_i$.
It follows that 
\begin{eqnarray*}
\|x^n-Tx^n\|_1&=&p_n^0+|p_n^1-p_n^0|+|p_n^2-p_n^1|+\cdots+|p_n^n-p_n^{n-1}|+p_n^n\\
&=&2\;\max\{p_n^i:0\leq i\leq n\}.
\end{eqnarray*}
Now, consider $n=2m$ Bernoullis trials, half of them with
success probability $\alpha_i={u\over m}$ and the other half
with $\alpha_i=1-{u\over m}$. Then 
$$\max\{p_n^i:0\leq i\leq n\}\geq p_{2m}^m=\PP(X=Y)$$
with $X,Y$ independent Binomials $B(m,{u\over m})$. When 
$m\to\infty$ these Binomials converge to Poissons
so that $p_{2m}^m$ tends to $\sum_{k=0}^\infty ({\exp(-u)u^k\over k!})^2=\exp(-2u)I_0(2u)$.
Since $\sqrt{\sum_{i=1}^{2m}\alpha_i(1\!-\!\alpha_i)}$ tends to $\sqrt{2u}$, it follows that
$p^m_{2m}\sqrt{\sum_{i=1}^{2m}\alpha_i(1\!-\!\alpha_i)}$
can be made as close as desired to the value $\eta=\max_{x\geq 0}\sqrt{x}\exp(-x)I_0(x)$.
Hence the optimal $\kappa$ lies in the interval $[\eta,\frac{1}{\sqrt{\pi}}]\sim[0.4688,0.5642]$
which leaves a margin of at most 17\%.

\section{Two direct applications of Theorem \ref{main}}

\subsection{Convergence of the iterates}
The following result, which is basically known ({\em cf.} \cite{bgk1,bgk2,gro,ish,bgk3,rei}),
shows how Theorem \ref{main} can be used to obtain the convergence 
of the iterates, proving at the same time the existence of fixed points.

\begin{proposition} Suppose $\sum\alpha_k(1\!-\!\alpha_k)=\infty$ and $x_k$ bounded. \\[0.5ex]
(a) If $x_k$ is relatively compact then $x_k\to\bar x$ for some $\bar x\in{\rm Fix}(T)$.\\[0.5ex]
(b) If $X$ is a Hilbert space then $x_k\rightharpoonup\bar x$ for some $\bar x\in{\rm Fix}(T)$.
\end{proposition}
\proof (a) Choose a convergent subsequence $x_{k_n}\to\bar x$. From (\ref{bnd0}) we obtain
$x_k-Tx_k\to 0$ so that $\bar x$ must be a fixed point. Since 
$$\|x_k\!-\!\bar x\|=\|(1\!-\!\alpha_k)(x_{k-1}\!-\!\bar x)+\alpha_k(Tx_{k-1}\!-\!T\bar x)\|\leq\|x_{k-1}\!-\!\bar x\|$$
we conclude that $\|x_k\!-\!\bar x\|$ decreases to 0.
\vspace{0.5ex}

(b) Since $I-T$ is maximal monotone and $x_k-Tx_k\to 0$, all weak cluster points of $x_k$ belong 
to ${\rm Fix}(T)$. As before $\|x_k-\bar x\|$ 
converges for all $\bar x\in{\rm Fix}(T)$ so that weak convergence follows from Opial's lemma.
\endproof

 \subsection{Unbounded domains}\label{stres}

 When $C$ is unbounded (\ref{bnd}) says nothing. However, if ${\rm Fix}(T)\neq\phi$ is 
 nonempty\,\footnote{A necessary and 
 sufficient condition to have ${\rm Fix}(T)\neq\phi$ is that the iterate sequence $\{x_k\}$ remains bounded
 ({\em cf.} \cite{rei0}).}, then for each
 $y\in \mbox{Fix}(T)$ we may still apply (\ref{bnd}) on the bounded subset $\tilde C=C\cap B(y,\|y\!-\!x_0\|)$ which
 satisfies $T(\tilde C)\subseteq \tilde C$ and $\mbox{diam}(\tilde C)\leq 2\|y-x_0\|$.
Hence, setting $\tilde\kappa=2\kappa$ and taking the infimum over $y\in\mbox{Fix}(T)$ we obtain
 \begin{equation}\label{bnd3}
\|x_n-Tx_n\|\leq \tilde\kappa\frac{\mbox{dist}(x_0,\mbox{Fix}(T))}{\sqrt{\sum_{i=1}^n\alpha_i(1\!-\!\alpha_i)}}.
\end{equation}

In particular, Theorem \ref{main} implies that (\ref{bnd3}) holds with $\tilde\kappa=\mbox{\scriptsize $\frac{2}{\sqrt{\pi}}$}\sim 1.1284$. 
In Hilbert spaces, \cite[Vaisman]{vai} established a sharper bound with 
$\tilde\kappa=1$.  We present this result which exploits
the well-known identity 
\begin{equation}\label{hi}
\|(1\!-\!\alpha)u+\alpha v\|^2=(1\!-\!\alpha)\|u\|^2+\alpha\|v\|^2-\alpha(1\!-\!\alpha)\|u-v\|^2.
\end{equation}

\begin{proposition} Let $T:C\to C$ be non-expansive on a convex $C\subset E$ with $E$ a Hilbert space
and $\mbox{\rm Fix}(T)$ nonempty. Then {\em (\ref{bnd3})} holds with $\tilde\kappa=1$. 
\end{proposition}
\proof It is known that $\|x_k-Tx_k\|$ decreases with $k$. Indeed,
\begin{eqnarray*}
\|x_k-Tx_k\|&=&\|(1\!-\!\alpha_k)x_{k-1}+\alpha_kTx_{k-1}-Tx_k\|\\
&\leq&(1\!-\!\alpha_k)\|x_{k-1}-Tx_{k-1}\|+\|Tx_{k-1}-Tx_k\|\\
&\leq&(1\!-\!\alpha_k)\|x_{k-1}-Tx_{k-1}\|+\|x_{k-1}-x_k\|\\
&=&(1\!-\!\alpha_k)\|x_{k-1}-Tx_{k-1}\|+\alpha_k\|x_{k-1}-Tx_{k-1}\|\\
&=&\|x_{k-1}-Tx_{k-1}\|.
\end{eqnarray*}
Now, using (\ref{hi}), for each $y\in\mbox{Fix}(T)$ we get
\begin{eqnarray*}
\|x_i\!-\!y\|^2\!\!&=&\!\|(1\!-\!\alpha_i)(x_{i-1}\!\!-\!y)+\alpha_i(Tx_{i-1}\!\!-\!Ty)\|^2\\
&=&\!\!(1\!-\!\alpha_i)\|x_{i-1}\!\!-\!y\|^2\!\!+\!\alpha_i\|Tx_{i-1}\!\!-\!Ty\|^2\!\!-\!\alpha_i(1\!-\!\alpha_i)\|x_{i-1}\!\!-\!Tx_{i-1}\|^2\\
&\leq&\!\|x_{i-1}-y\|^2-\alpha_i(1\!-\!\alpha_i)\|x_{i-1}-Tx_{i-1}\|^2.
\end{eqnarray*}
Summing these inequalities we see that 
$$\mbox{$\sum_{i=1}^n$}\alpha_i(1\!-\!\alpha_i)\|x_{i-1}-Tx_{i-1}\|^2\leq  \|x_0-y\|^2-\|x_n\!-\!y\|^2$$
and the monotonicity of  $\|x_k-Tx_k\|$ yields
$$\|x_n-Tx_n\|\sqrt{\mbox{$\sum_{i=1}^n\alpha_i(1\!-\!\alpha_i)$}}\leq  \|x_0-y\|.$$
The conclusion follows by taking the infimum over $y\in\mbox{Fix}(T)$.
\endproof

\vspace{1ex}\noindent
{\sc Remark:} The previous proof yields a slightly sharper estimate
$$\|x_{n-1}-Tx_{n-1}\|\leq \frac{\mbox{dist}(x_0,\mbox{Fix}(T))}{\sqrt{\sum_{i=1}^n\alpha_i(1-\alpha_i)}}$$
with $x_{n-1}$ in place of $x_n$ on the left.

\section{Auxiliary results}\label{auxiliary}
\subsection{A Hoeffding-type inequality}\label{hoeffding}
In this short section we establish a Hoeffding-type inequality for sums of Bernoullis
and Poisson variables. We consider an integer function $g:\NN\to\RR$
satisfying the convexity inequalities $g(k)\leq\frac{1}{2}[g(k\!-\!1)+g(k\!+\!1)]$ for all $k\geq 1$.
\begin{proposition}\label{poisson}
Let  $S=X_1+\cdots+X_m$ be a sum of independent Bernoulli trials
with success probabilities $\PP(X_i\!=\!1)\!=\!p_i$, and let $z=\EE(S)=p_1+\ldots+p_n$.
Then $\EE[g(S)]\leq\EE[g(Z)]$ where $Z\sim P(z)$
is a Poisson with the same mean.
\end{proposition}
\proof Let us first note that the expected value $\EE[g(S)]$
increases if we replace any variable $X_i$ by a sum $X_i'+X_i''$ 
of independent Bernoullis with 
$$\PP(X_i'=1)=\PP(X_i''=1)=\mbox{$\frac{p_i}{2}$}.$$
Indeed, for $k\!\in\!\NN$ let $A(k)\!=\!\EE[g(k\!+\!X_i)]$ and $B(k)\!=\!\EE[g(k\!+\!X_i'\!+\!X_i'')]$
so that
\begin{eqnarray*}
A(k)&=&(1-p_i) g(k)+p_i g(k+1)\\
B(k)&=&\mbox{$(1-\frac{p_i}{2})^2 g(k)+p_i(1-\frac{p_i}{2}) g(k+1)+(\frac{p_i}{2})^2 g(k+2)$}.
\end{eqnarray*}
Taking their difference we have
$$B(k)-A(k)=\mbox{$(\frac{p_i}{2})^2[g(k)-2 g(k+1)+ g(k+2)]$}\geq 0$$
so that replacing $k$ by the random variable $\sum_{j\neq i}X_j$
and taking expectation we obtain the asserted monotonicity.

Now, a well-known result by Hoeffding \cite[Theorem 3]{hoef} 
proves that\footnote{As a matter of fact, Hoeffding assumes $g$ strictly convex but the general case 
follows by applying his result to $g(x)+\epsilon x^2$ with $\epsilon\downarrow 0$.} 
$\EE[g(S)]\leq\EE[g(S_1)]$ with $S_1\sim B(n,p)$ a binomial
with $p=\frac{1}{n}(p_1+\ldots+p_n)$.
Writing $S_1$ as a sum of $n$ Bernoullis $B(p)$
and sequentially replacing each term  by two Bernoullis 
$B(\frac{p}{2})$, the expected value increases in each step
and we get $\EE[g(S)]\leq\EE[g(S_2)]$ with $S_2\sim B(2n,p/2)$.
Iterating this doubling argument we obtain $\EE[g(S)]\leq \EE[g(S_k)]$ 
where $S_k\sim B(2^kn,p/2^k)$. Since $\EE(S_k)=z$ for all $k$, the 
result follows by letting $k\to\infty$ and noting that $S_k$ converges to a Poisson variable 
$Z\sim P(z)$.
\endproof

\subsection{An identity for Catalan numbers}\label{acatalan}
In proving Proposition \ref{equal} we used the identity 
$$C_k=\mbox{$\sum_{j=0}^k(-1)^j2^{k-j}{k\choose j}{j\choose \lfloor j/2\rfloor}$}.$$
Since this is not found in standard textbooks, for completeness we provide a proof. 
For each $a \in\ZZ$ and $P(x)$ a Laurent polynomial  ({\em i.e.} a function whose Laurent series has 
finitely many terms) we denote by $[x^a]P(x)$ the coefficient of $x^a$ in $P(x)$. 
We observe that for each non-negative integer $j$ we have
\begin{align*}
[x^0](x^2\!+\! x^{-2})^j&= 
\left\{\begin{array}{cl} 
~\binom{j}{\frac{j}{2}}&\text{~ for $j$ even}\\[1.5ex]
~0&\text{~ for $j$ odd}\end{array} \right.\\
[x^2](x^2\!+\! x^{-2})^j&= 
\left\{\begin{array}{cl} 
0&\text{ for $j$ even}\\[0.5ex]  \binom{j}{\frac{j-1}{2}}&\text{ for $j$ odd}
\end{array} \right.
\end{align*}
so we can write $\binom{j}{\lfloor j/2 \rfloor}=([x^0]+[x^2]) (x^2 + x^{-2})^j$ and therefore
\begin{align*}
\mbox{$\sum_{j=0}^k(-1)^j2^{k-j}{k\choose j}{j\choose \lfloor j/2\rfloor}$}&= ([x^0]+[x^2])\, \mbox{$\sum_{j=0}^k$}\mbox{$\binom{k}{j}2^{k-j} $}(- x^2\!-\!x^{-2})^j\\
&= ([x^0]+[x^2])\, \left(2\! -\! x^2\!-\!x^{-2}\right)^k \\
&=([x^0]+[x^2])\, \left(-(x^1\!-\!x^{-1})^2\right)^k \\
&= ([x^0]+[x^2])\, (-1)^k\left(x^1\! -\! x^{-1}\right)^{2k}\\
&=\mbox{$ \binom{2k}{k}-\binom{2k}{k+1}$} = C_k.
\end{align*}

\end{document}